\documentclass[journal]{IEEEtran}
\usepackage{multirow}
\usepackage{multicol}
\usepackage{lipsum}
\usepackage{graphicx}
\usepackage{graphics}
\usepackage{amsmath}
\usepackage{amsbsy}
\usepackage{blindtext}
\usepackage{tcolorbox}
\usepackage{amssymb}
\usepackage{scalerel}
\usepackage{mathabx}
\usepackage{verbatim}
\usepackage{booktabs}
\usepackage{rotating,tabularx}
\usepackage{mathrsfs} 
\usepackage{hyperref}

\usepackage{adjustbox}
\usepackage{soul}
\usepackage{xcolor}
\usepackage{cite}
\usepackage{tikz}
\usepackage{color, colortbl}
\usepackage[first=0,last=9]{lcg}
\definecolor{LightGray}{rgb}{0.7,0.7,0.7}

\usepackage[wby]{callouts}
\usetikzlibrary{shapes.multipart}

\usepackage{array}
\usepackage{makecell}

\allowdisplaybreaks

\usepackage{amsthm}
\theoremstyle{definition}

\theoremstyle{remark}

\usepackage[utf8]{inputenc}
\usepackage[english]{babel}

\usepackage[caption=false,font=footnotesize]{subfig}
\usepackage[T1]{fontenc}
\usepackage{scalerel,stackengine}
\newcommand\reallywidecheck[1]{%
\savestack{\tmpbox}{\stretchto{%
  \scaleto{%
    \scalerel*[\widthof{\ensuremath{#1}}]{\kern-.6pt\bigwedge\kern-.6pt}%
    {\rule[-\textheight/2]{1ex}{\textheight}}%WIDTH-LIMITED BIG WEDGE
  }{\textheight}% 
}{0.5ex}}%
\stackon[1pt]{#1}{\scalebox{-1}{\tmpbox}}%
}
\stackMath
\IEEEoverridecommandlockouts
\IEEEoverridecommandlockouts

\renewcommand{\baselinestretch}{1}

\newcommand*{\mrn}{\textcolor{black}}
\newcommand*{\rn}{\textcolor{black}}

\newif\ifarxiv
\arxivtrue
\arxivfalse

\begin{document}

\title{\LARGE\bf
\rn{Designing Sparse AC False Data Injection Attack}}

\author{Mohammadreza Iranpour$^{\ast}$, Mohammad Rasoul Narimani$^{\ast}$
\thanks{${\ast}$: Department of Electrical and Computer Engineering, California State University Northridge (CSUN). Rasoul.narimani@csun.edu. Support from NSF contract \#2308498.}% 
}

\maketitle

\begin{abstract}

\mrn{False Data Injection (FDI) attacks pose significant threats by manipulating measurement data, leading to incorrect state estimation. Although numerous studies have focused on designing DC FDI attacks, few have addressed AC FDI attacks due to the complexity of incorporating non-linear AC power flows in the design process. Additionally, designing a sparse AC FDI attack presents another challenge because it involves solving a mixed-integer nonlinear programming problem with nonconvex constraints, which is inherently difficult. This paper explores the design and implementation of a sparse AC FDI attack, where the attacker strategically selects a minimal set of measurements to manipulate while maintaining the nonlinearity and interdependence of AC power flow equations. The objective is to minimize the number of altered measurements, thereby reducing the attack's detectability while achieving the desired state estimation error. The problem is formulated as a Mixed Integer Nonlinear Programming (MINLP) problem. Binary variables indicate the selection of measurements to be manipulated, and continuous variables represent the measurement values. 
%Power flow constraints dynamically change depending on whether the measurements are selected or remain fixed. 
An optimization problem is designed to minimize the number of binary variables, translating into a sparse attack, while ensuring the attack remains efficient and hard to detect. The big-M method and conditional constraints are utilized to handle the fixed and variable measurement parameters effectively. Simulation results on the standard IEEE 57-bus test system demonstrate the efficacy of the sparse AC FDI attack in terms of its impact on state estimation and the minimal number of measurements required for successful implementation.}

\end{abstract}

\section{Introduction}
\label{Introduction}
\mrn{The rapid digitalization and automation of power systems have significantly improved their efficiency and reliability. However, this increased dependence on digital infrastructure also exposes the power grid to potential cyber-threats, such as False Data Injection (FDI) attacks. These attacks aim to disrupt the normal operation of power systems by manipulating data reported by sensors and meters~\cite{liu2011false}. When the attack vector meets specific conditions, the adversary can introduce errors in state estimation, bypassing standard residue-based Bad Data Detection (BDD) tests. This leads to incorrect state estimation and potentially operational and economic impacts~\cite{du2021targeted, boyaci2021joint,boyaci2022generating, boyaci2022infinite, boyaci2022spatio}.} 

\mrn{Traditionally, research on FDI attacks has predominantly focused on DC FDI models due to their linear nature, which simplifies the attack design and analysis process. However, real-world power systems operate under AC conditions, characterized by nonlinear power flow equations, making AC FDI attacks a more accurate and practical area of study. Despite their relevance, the complexity of incorporating nonlinear AC power flows into attack design has resulted in limited research on AC FDI attacks~\cite{rahman2013false}. However, numerous researchers have extensively studied FDI attacks and explored their background, construction methods, detection, and defense strategies. Researchers have looked at different ways to create FDI attacks, including those with limited resources, those based on incomplete system knowledge, and data-driven methods~\cite{zhang2019false}. Attacks with limited budgets focus on using minimal effort, like compromising a small number of meters, known as the k-sparse problem. Systematic methods have shown that FDI attacks can be created by solving an optimization problem to find a stealthy and sparse data injection vector that leads to false state estimates~\cite{zhao2018generalized,zhang2018can,du2021targeted}.}

\mrn{From an optimized FDI attack perspective, numerous studies have focused on optimization approaches that emphasize sparsity and minimizing the number of compromised measurements in power systems. Some researchers have used numerical and heuristic methods, while others have employed analytical methods. A reduced row-echelon (RRE) form-based greedy method is used to compute the minimum number of targeted measurements needed for an attack~\cite{nayak2020modelling}. Additionally, genetic algorithms and neural networks have been applied to construct the least-effort attack vector~\cite{jin2017semidefinite}.
The minimum number of sensors required for stealth FDI attacks can be quantified by formulating a minimum cardinality problem, with various algorithms proposed for efficient computation. These problems, studied extensively in the past decade, are characterized by nonconvexity and NP-hardness. Different methods have been presented to solve them. For instance, a two-step approach for computing sparse attack vectors based on successive, iteratively reweighted, L1-norm convex relaxations is proposed in~\cite{kim2011strategic}. In~\cite{ozay2013sparse}, sparse attack vectors are created by incorporating a regularizing L1-norm term in the objective function or explicitly including a cardinality condition for the attack vector in the constraints.}
\mrn{When the sparse subset of measurements is devised based on a prior specification of a non-zero entry in the attack vector, methods such as mixed-integer linear programming (MILP), L1-norm convex relaxation, or a min-cut scheme can be employed~\cite{teixeira2015secure}. In~\cite{alexopoulos2020complementarity}, both exact and relaxed complementarity reformulations of cardinality minimization are proposed to identify the minimal sets of measurements that need to be maliciously modified for a successful FDIA by solving a mixed-integer program with half-complementarity constraints~\cite{feng2013complementarity}. Other methods, such as branch-and-cut~\cite{de2001branch, bennett1997parametric}, can also be used to design sparse FDI attacks.}

\mrn{Although these studies have analyzed sparse DC FDI attacks, there is limited research on designing sparse AC FDI attacks. Designing a sparse AC FDI attack is a challenging optimization problem because the nonlinear power flow equations governing the physics of the power system must be incorporated into the optimization process\cite{narimani2023tightening,narimani2020tightening, narimani2020strengthening,narimani2018comparison, narimani2018empirical,narimani2018improving}. A sparse attack aims to alter the minimum number of measurements needed to achieve the desired impact, thereby reducing the likelihood of detection. This involves solving a mixed-integer nonlinear programming (MINLP) problem with non-convex constraints, due to the nonlinear power flow equations in AC systems. Additionally, the attack design must account for selecting a minimal subset of measurements with variable parameters while keeping unselected measurements fixed.}

\mrn{Although \cite{hijazi2017convex} provides a well-designed MINLP framework for various power system problems, there has been no implementation of an MINLP formulation for sparse AC FDI attacks. To address these complexities, this paper introduces a process for designing an AC FDI attack based on nonlinear power flow equations. These equations are enforced based on their corresponding binary variables, dynamically adjusting depending on whether the measurements are selected or fixed, which facilitates an efficient and robust solution to sparse AC FDI attack problems.}
\mrn{In this connection, we first define the process of designing an AC-based FDI attack that satisfies AC power flow equations in Section \ref{sec:AC design}. We then formulate an optimization problem to design a sparse AC FDI attack that both satisfies AC power flow equations and selects a minimum set of measurements for manipulation, as detailed in Section \ref{sec:sparse AC design}. In Section \ref{sec:results}, we assess the capability of the proposed method by applying it to IEEE test systems. Section \ref{sec:conclusion} concludes the paper.}

\section{Process of designing the Sparse AC FDI attack}
\label{sec:AC design}

\mrn{Engineers use monitoring networks to keep power systems running smoothly by providing measurement data to state estimation platforms in control centers. State estimators analyze parameters like voltage levels and angles, which are essential for making decisions. Usually, accurate sensor readings lead to state estimates that closely match true values, but abnormal measurements can distort these estimates. To find inconsistencies, researchers calculate the measurement residual $r = z - h(x)$, which shows the difference between actual measurements ($z$) and estimated ones ($h(x)$). If $\|z - h(x)\|>\tau$, where $\tau$ is a set threshold, it signals erroneous measurements. This process is known as Bad Data Detection (BDD) in power systems.}

\mrn{This weakness in power transmission monitoring systems creates an opportunity for attackers to exploit vulnerabilities in the monitoring networks. Although bad data detectors are designed to identify and remove incorrect measurements, adversaries can inject carefully crafted false measurements that disrupt data distribution while avoiding detection. This serves as a starting point for designing a successful FDI attack.} 
\mrn{In this context, if a false variable state vector ($x_{attack}$) is injected into the measurements selected for conducting the attack, the observed measurements become $z_{attack} = z + a$, where $a$ is the attack vector. Consequently, the residual of the state-estimated variables under attack is $r_{attack} = z_{attack} - h(x_{attack})$. By substituting $z_{attack} = z + a$ and adding and subtracting $h(x)$, we can write $r_{attack} = z + a - h(x_{attack}) + h(x) - h(x)$. Since $r = z - h(x)$, we have $r_{attack} = r + a - h(x_{attack}) + h(x)$. If we define $a = h(x_{attack}) - h(x)$, then $r_{attack}$ will equal $r$. Thus, if the attack vector $a$ is defined as in Equation~\eqref{eq:vector a}, the residue-based bad data detection tests cannot detect the attack vector $a$ because the injected false data does not affect the residual.}

\vspace{-.3cm}
\mrn{\begin{align}
  a = \| h(x_{attack}) - h(x) \|.
\label{eq:vector a}
\end{align}}

\mrn{The main difference between designing FDI attacks for AC and DC systems lies in the function $h(x)$. In DC FDI attacks, this function uses constant coefficients from DC power flow equations. In contrast, AC FDI attacks use $h$ as a set of nonlinear power flow equations that relate measurements to state variables. To design a successful AC FDI attack, it's crucial to consider the nonlinear relationships in power flow equations that connect state variables with quantities like power flow in lines and power injections into buses for the selected measurements.}

\mrn{To determine the function $h(x_{attack})$, several key assumptions must be considered. First, any changes in power flow equations resulting from the injected false data need to be rationalized to maintain consistency. It is important to clearly distinguish the attack area, where modifications occur, from the unaffected normal area. The attack area should be a contiguous set of buses, while the normal area remains unchanged. Designate primary focal buses to define the attack zone and expand it by including zero-injection buses and neighboring non-zero injection buses as boundaries. Additionally, the total power transfers between the affected and unaffected regions must be kept equivalent during an FDI attack. Surround the attack area with buses that have power injections to manage power changes within the zone, and ensure that state variables at the boundary of the attack area remain unchanged to confine all modifications within the attack zone. Finally, the algebraic sum of generated and consumed power within the attack zone must be preserved. For zero-injection buses, the sum of active and reactive power flows should be zero, while for non-zero injection buses, the post-attack power injection should equal the primary injection power plus the sum of power flow changes in connected lines within the attack zone. These assumptions are essential for designing AC FDI attacks.}

\mrn{To design a sparse AC FDI attack, we need to introduce an additional discrete variable to select the minimum set of allowable measurements for a successful attack. In the process of sparse FDI attacks, we work with PMU measurements, including voltage magnitude and angle measurements. In this context, unselected measurements within the predefined attack zone should be fixed at their corresponding values before the attack is implemented, while the state of selected measurements (voltage magnitude and angle) should remain variable during the solving process. Based on these descriptions, we formulate the appropriate functions $h(x)$ that ensure a successful AC FDI attack while considering the minimum set of measurements to be manipulated. This formulation is detailed in the next section. }

\section{sparse AC FDI attack formulation}
\label{sec:sparse AC design}

\mrn{This section formulates the problem of selecting a sparse set of measurements for manipulation to conduct an AC FDI attack in power system, as presented in the section~\ref{sec:AC design}.}
\mrn{The key challenge in designing a sparse AC FDI attack lies in managing both continuous and binary variables, as well as dealing with nonlinear equations. This transforms the problem into a mixed-integer nonlinear programming (MINLP) problem, which is one of the most complex types of optimization problems. To formulate the problem, we define a binary variable $z$ to determine which PMU measurements should be manipulated and use the big M method\cite{hijazi2017convex} to handle constraints involving binary variables.
In fact, by selecting a sparse set of measurements on the buses of the system, we can adjust the voltage magnitude and angle at these buses to meet the assumptions necessary for implementing a successful attack scenario, as described in the section~\ref{sec:AC design}. The goal is to design a sparse FDI attack, and binary variables have been defined to determine the smallest set of PMU measurements for manipulation.
%The second point is considering all power flow and transferring assumptions in formulating this problem which will be considered as constraints of the objective function of problem. 
Continuing, we modeled these assumptions as an optimization problem to calculate the attack vector $a$. Suppose the sets of buses and lines of a system are represented by $\mathcal{B}$, and $\mathcal{L}$, respectively. Also, $\mathcal{B_A}$, and $\mathcal{L_A}$, are corresponding sets of buses and lines in the attack zone,  respectively. Let $S_m = {P}_m + j {Q}_m$ represents the complex power injection, $V_m$ and $\theta_m$ represent the voltage magnitude and angle at bus~$m\in\mathcal{B_A}$, each line $\left(m,l\right)\in\mathcal{L_A}$ is modeled as a $\Pi$ circuit with mutual admittance $g_{ml}+j b_{ml}$ and shunt admittance $j b_{c,ml}$ and the voltage angle difference between buses $m$ and $l$ for $(m,l)\in\mathcal{L}$ is denoted as $\theta_{ml}=\theta_{m}-\theta_{l}$.}

%The difference between state variables, including both voltage magnitudes and angles and their values prior to conducting FDI attack is shown in Equation~\eqref{eq:optimization} as vector $c$.

%\mir{\begin{equation}
%\label{eq:optimization}
% c= [\tilde{V}_{m}- V_{m,fix}, \tilde{\theta}_{m}-\theta_{m,fix}]
%\end{equation}}

\mrn{Also, consider $\tilde{V}_{m}$ and $\tilde{\theta}_{m}$ aas variables representing the voltage magnitude and angle that need to be calculated to conduct the attack, while $V_{m,fix}$ and $\theta_{m,fix}$ are the known values of voltage magnitude and angle before the attack. We aim to fix the voltage magnitude and angle of unselected measurements at $V_{m,fix}$ and $\theta_{m,fix}$. Meanwhile, we want to minimize the number of selected measurements for manipulation. By minimizing the sum of selected measurements, we can define an objective function that considers constraints to satisfy all assumptions presented in Sec \ref{sec:AC design}, as represented in Equations  \eqref{eq:obj1}-\eqref{eq:reactive_overload}.}

\mrn{\begin{small}
\begin{subequations}
\begin{align}
%&\min\quad \textstyle\sum_{\small{{m}\in \mathcal{B_A}}}
% \left(\tilde{V}_{m}- V_{m,fix})^2 + (\tilde{\theta}_{m}-\theta_{m,fix})^2 \label{eq:obj} \right)\\
&\min\quad \textstyle\sum z_{i}
 \label{eq:obj1}\\
%&\quad\qquad\qquad+c_{1,k}\left( {\tilde P_k^g\cos (\psi_l ) - \tilde Q_k^g\sin (\psi_l )} \right)+c_{0,k}\\
%&\nonumber \min\quad \textstyle\sum_{\footnotesize{{k}\in \mathcal{G}}}
%c_{2,k}\left( {\tilde P_k^g\cos (\psi_l ) - \tilde Q_k^g\sin (\psi_l)} \label{eq:RQC obj} \right)^2\\
%&\quad\qquad\qquad+c_{1,k}\left( {\tilde P_k^g\cos (\psi_l ) - \tilde Q_k^g\sin (\psi_l )} \right)+c_{0,k}\\
&\nonumber \text{subject to} \quad \left(\forall i\in\mathcal{B_A}, \forall   \left(l,m\right) \in\mathcal{L_A}\right)\\
&V \geq V_{l}Z+ (1-Z)V_{m,fix}\label{eq:Zg_VOLTAGE} \\
&V \leq V_{u}Z + (1-Z)V_{m,fix}\label{eq:Zl_VOLTAGE}\\
&\theta\geq -MZ+ (1-Z)\theta_{m,fix}\label{eq:Zg_angle}\\ &\theta\leq MZ+ (1-Z)\theta_{m,fix}\label{eq:Zl_angle}\\
&-\pi /6\leqslant  \theta_{lm}\leqslant \pi /6\label{eq:Z_difangle}\\
&\!\!\!\!\!\!\! g_{sh,i}\, \tilde{V}_i^2+\sum_{\substack{(l,m)\in \mathcal{L},\\\text{s.t.} \hspace{3pt} l=i}} \!\tilde{P}_{lm}+\!\!\sum_{\substack{(l,m)\in \mathcal{L},\\\text{s.t.} 
\label{eq:active_injection}\hspace{3pt} m=i}} \!\!\tilde{P}_{ml}= P_{i,G}-P_{i,D}, \\
&\!\!\!\!\!\!\! -b_{sh,i}\, \tilde{V}_i^2+\!\!\!\!\!\!\sum_{\substack{(l,m)\in \mathcal{L},\\ \text{s.t.} \hspace{3pt} l=i}} \!\!\tilde{Q}_{lm}+\!\!\!\sum_{\substack{(l,m)\in \mathcal{L},\\ \text{s.t.}\label{eq:reactive_injection}\hspace{3pt} m=i}} \!\!\!\!\tilde{Q}_{ml}=Q_{i,G}-Q_{i,D},\\
&\nonumber\!\!\!\!\!\!\! \tilde{P}_{lm} \!=\! g_{lm} \tilde{V}_l^2\! -\! g_{lm} \tilde{V}_l V_m\cos\left(\tilde{\theta}_{l}-\theta_{m}\right)\!\\
&\qquad\qquad -\! b_{lm} \tilde{V}_l V_m\sin\left(\tilde{\theta}_{l}-\theta_{m}\right),\\
\label{eq:qik1}
&\!\!\!\!\!\!\! \nonumber \tilde{Q}_{lm} = -\left(b_{lm}+b_{c,lm}/2\right) \tilde{V}_l^2 + b_{lm} \tilde{V}_l V_m\cos\left(\tilde{\theta}_{l}-\theta_{m}\right)\\ &\qquad\qquad  - g_{lm} \tilde{V}_l V_m\sin\left(\tilde{\theta}_{l}-\theta_{m}\right),\\
\label{eq:pki1}
&\nonumber \!\!\!\!\!\!\!\tilde{P}_{ml}\! =\! g_{lm} V_m^2\! -\! g_{lm} \tilde{V}_l V_m\cos\left(\tilde{\theta}_{l}-\theta_{m}\right)\!\\
&\qquad\qquad +\! b_{lm} \tilde{V}_l V_m\sin\left(\tilde{\theta}_{l}-\theta_{m}\right),\\
\label{eq:qki1}
&\!\!\!\!\!\!\!\nonumber \tilde{Q}_{ml} = -\left(b_{lm}+b_{c,lm}/2\right) V_m^2 + b_{lm} \tilde{V}_l V_m\cos\left(\tilde{\theta}_{l}-\theta_{m}\right)\\ &\qquad\qquad  + g_{lm} \tilde{V}_l V_m\sin\left(\tilde{\theta}_{l}-\theta_{m}\right)\\
\label{eq:active_overload}
&\!\!\!\!\!\!\! \tilde{P}_{ml} = W_{lm}*{P}_{lm}^{P
F},\\
\label{eq:reactive_overload}
&\!\!\!\!\!\!\! \tilde{Q}_{ml} = W_{lm}*{Q}_{lm}^{P
F}.
\end{align}
\end{subequations}
\end{small}}

\mrn{In equations \ref{eq:Zg_VOLTAGE} and \ref{eq:Zl_VOLTAGE}, when $z=0$, the voltage magnitude of the corresponding bus is fixed at its pre-attack value $V_{m,fix}$, and when $z=1$, the voltage magnitude of the corresponding bus is treated as a variable within the system's upper and lower voltage bounds. Similarly, in Equations~\ref{eq:Zg_angle} and \ref{eq:Zl_angle}, when $z=0$, the voltage angle of the corresponding bus is fixed at its pre-attack value $\theta_{m,fix}$, and when $z=1$, the voltage angle of the corresponding bus is considered a variable with unlimited bounds (since $M$ is set to a large number). In power system equations, power flows are related to voltage angle differences, and we have considered Equation \ref{eq:Z_difangle} to bound these values.}

\mrn{In these equations, 
$\tilde{P}_{ml}$, $\tilde{Q}_{ml}$, $\tilde{V}_{l}(\tilde{V}_{m})$ and $\tilde{\theta}_{l}(\tilde{\theta}_{m})$ represent the active and reactive power flow between buses 
$m$ and $l$ $(m,l)\in\mathcal{L}$ and the voltage values of buses $m$ and 
$l$ after the attack. Similarly, 
$P_{ml}$, $Q_{ml}$, $V_{l}(V_{m})$ and $\theta_{l}(\theta_{m})$ are the corresponding quantities before the attack. 
$P_{m,G}$, $Q_{m,G}$, $P_{m,D}$, and $Q_{m,D}$ denote the active and reactive power generation and the active and reactive power demand at bus $m$ , respectively. It is notable that in Equations~\eqref{eq:active_injection} and~\eqref{eq:reactive_injection}, when these equations are applied to zero-injection buses, the right-hand side of these equations equals zero.}
\mrn{Equations \ref{eq:active_overload} and \ref{eq:reactive_overload} are additional constraints for overloading a specified line in the attack zone with a predefined coefficient $W$. By considering this constraint, we can design an optimal AC FDI attack with a specific aim, such as overloading a particular line by a predefined coefficient. ${P}_{lm}^{P
F}$ and ${Q}_{lm}^{P
F}$ represent the active and reactive power flow before the attack. After solving this optimization problem, we can calculate the power injection values of non-zero injection buses within the attack zone as follows:}

\vspace{-.2cm}
\mrn{\begin{subequations}
\begin{small}
\begin{align}
\label{eq:non_zero_injection}
&\tilde{P}_{m}=P_{m}+\sum_{(m,l)\in \mathcal{L}_A}(\tilde{P}_{m,l}-P_{m,l}),\\
&\tilde{Q}_{m}=Q_{m}+\sum_{(m.l)\in \mathcal{L_A}}(\tilde{Q}_{m,l}-Q_{m,l}).
\end{align}
\end{small}
\end{subequations}}

\mrn{Here, $\tilde{P_{m}}$, $\tilde{Q_{m}}$, $P_{m}$, and $Q_{m}$ represent the active and reactive power injections at bus $m$ after and before the attack, respectively. Notably, in Equations~\eqref{eq:non_zero_injection}, $\tilde{P}{m,l}$ and $\tilde{Q}{m,l} \in \mathcal{L}_A$ should be considered only for the lines within the attack zone that are connected to bus $m$. After calculating all of these values, we can define the attack vector $a$ based on Equation~\eqref{eq:vector a} as follows:}

\mrn{\begin{small}
\begin{align}
\label{eq:attack_vector_final}
 &\nonumber a= [\tilde{P_{ml}}- P_{ml}, 
     \tilde{Q_{ml}}- Q_{ml},
     \tilde{P_{m}}- P_{m},\\
  &~~~~~~\nonumber 
 \tilde{Q_{m}}- Q_{m},
     \tilde{V_{m}}\angle \tilde{\theta_{m}}- V_{m}\angle\theta_{m}]^T
\end{align}
\vspace{-.3cm}
\end{small}}

\section{Assessing the proposed method and numerical results}
\label{sec:results}
\begin{figure}
    \centering
\includegraphics[scale=0.47,trim= 1.5cm 8cm 2cm 1cm,clip]{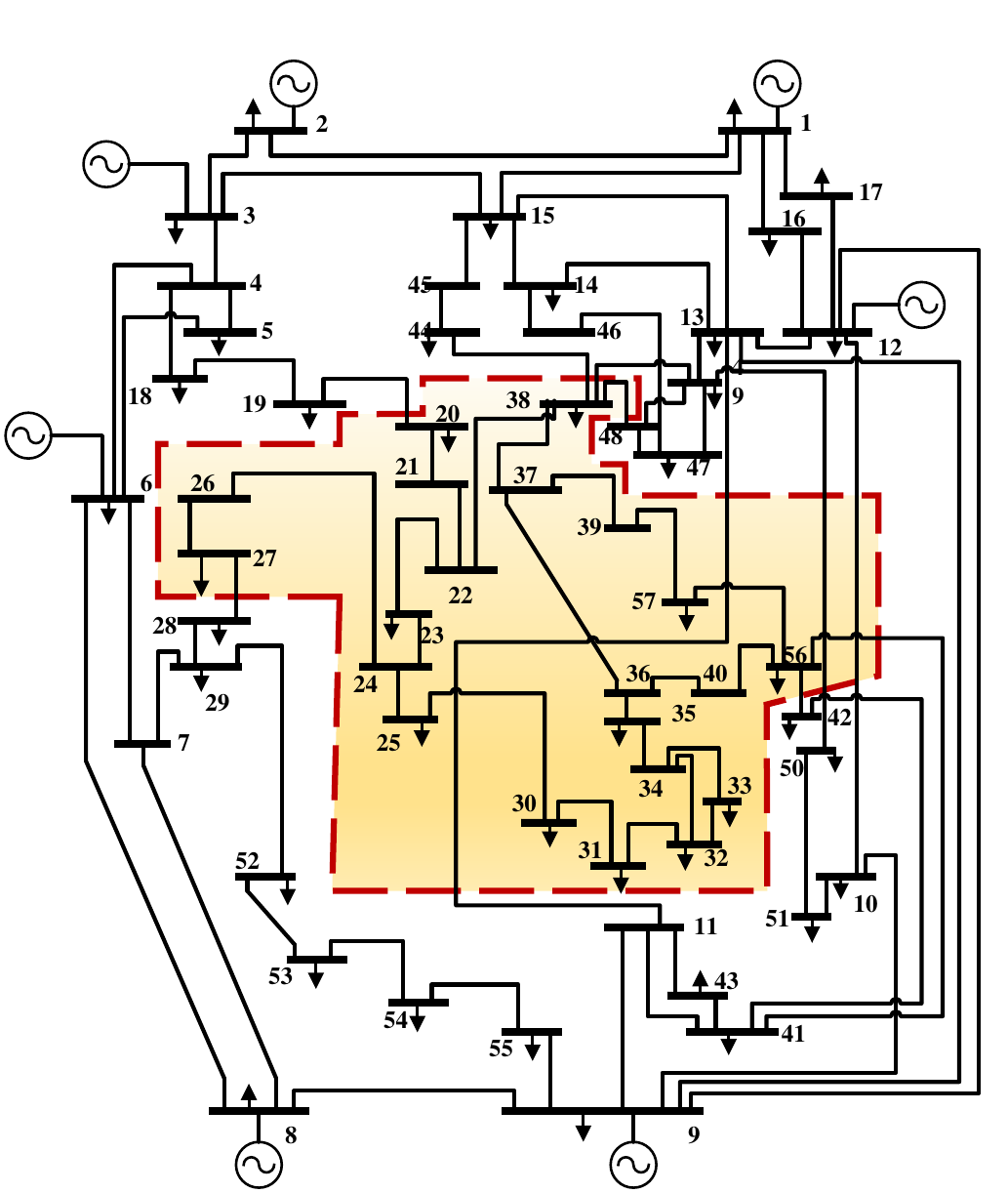}
	\caption{One-line diagram depicting the IEEE 57-bus test system, with the
attack zones indicated}
	\label{fig:attack_zone_57}
 \vspace{-.5cm}
\end{figure}

\mrn{In this section, we apply the proposed approach for designing sparse AC-FDI attacks to the IEEE 57-bus test system from the PGLib-OPF v18.08 benchmark library~\cite{pglib} to evaluate its effectiveness. We consider two scenarios: designing arbitrary and sparse AC-FDI attacks. In the arbitrary attack scenario, we demonstrate that, after identifying the attack zone, all buses except the boundary buses are selected for manipulation. This includes PMU measurements for voltage magnitude and angle at these buses. In contrast, the sparse attack scenario aims to identify a minimal set of measurements within the attack zone that can effectively alter the system. In this scenario, we determine the optimal number of measurements required to conduct a successful attack. The IEEE 57-bus test system and the selected attack zone for the sparse attack scenario are illustrated in Fig.~\ref{fig:attack_zone_57}. This visualization highlights the configuration used to implement the sparse attack and validate the proposed approach.}

\textbf{Arbitrary Attack:}
\mrn{In this scenario, we consider buses 20, 21, 22, 23, 24, 25, 26, 27, 30, 31, 32, 33, 34, 35, 36, 37, 38, 39, 40, 56 and $57$ as part of the attack zone. Within this zone, buses 
$20, 27, 38, 38,$ and 
$56$ are designated as boundary buses. For this setup, we designed the FDI attack by substituting a fixed quantity into the objective function of Equation~\eqref{eq:obj1} and enforcing the constraint in Equation~\ref{eq:active_overload} with $W=3$. This constraint ensures that the flow in the line between buses $23$ and $24$ is increased to three times its original flow before the attack.
In this scenario, there is no requirement to select a minimal set of measurements. Therefore, all buses in the attack zone, except for the boundary buses, were chosen to have their measurement values altered. As shown in Table \ref{Table:arbitrary}, the voltage magnitude and angle for all buses (excluding the boundary buses) within the attack zone have been modified. In this table, gray rows represent boundary buses with fixed voltage, while pink rows denote flexible nodes within the attack zone that were selected for manipulation after solving the problem.}

\begin{figure}
    \centering
\includegraphics[scale=0.47]{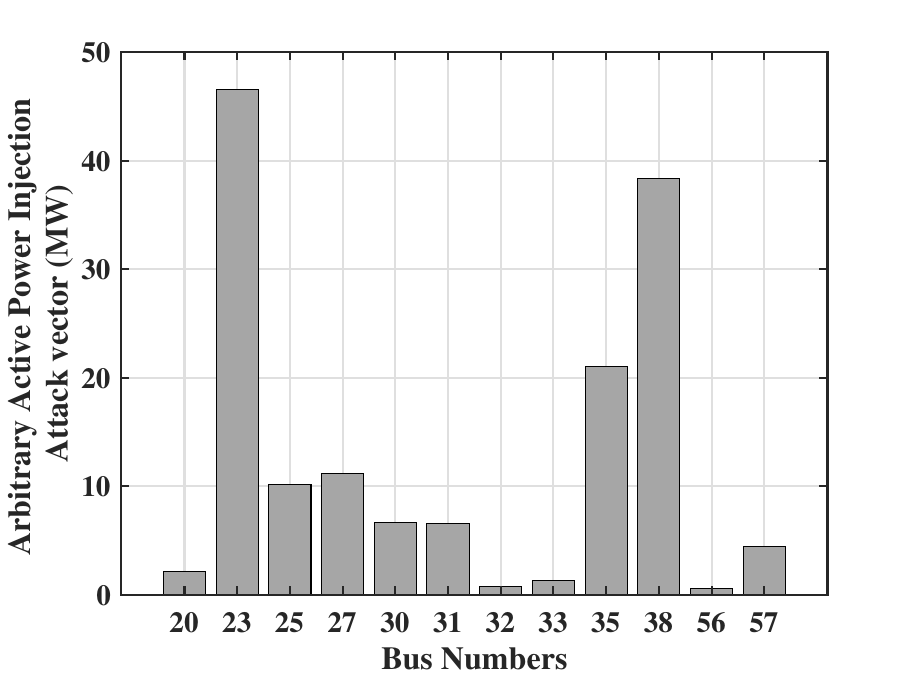}
	\caption{Active Power injection attack vector for arbitrary scenario.}
	\label{fig:powerinjection_arbitrary}
  %\vspace{-.4cm}
\end{figure}

\begin{figure}
 \vspace{-.4cm}
    \centering
\includegraphics[scale=0.47]{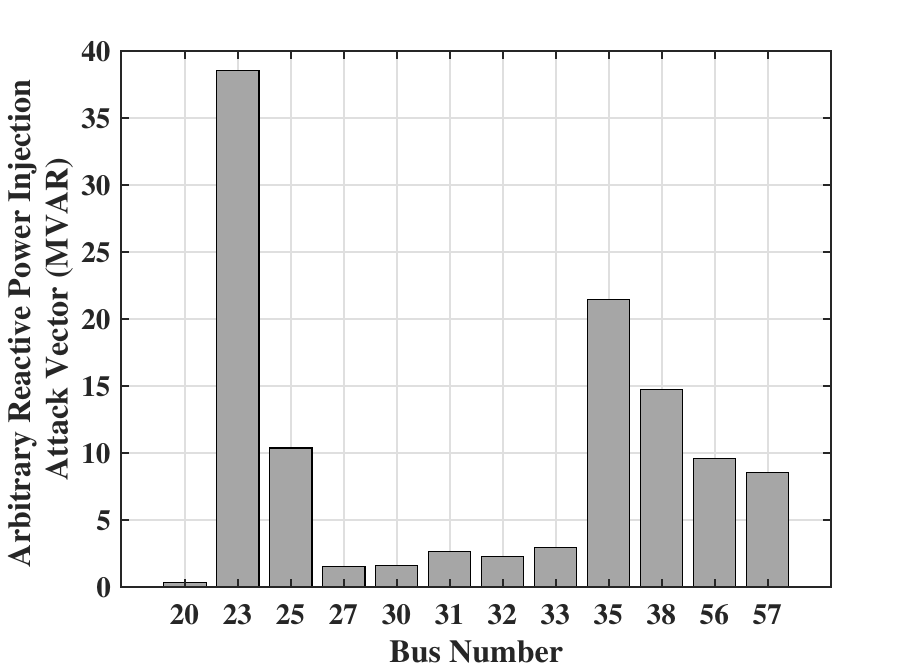}
	\caption{Reactive Power injection attack vector for arbitrary scenario.}
	\label{fig:reactivepowerinjection_arbitrary}
 %\vspace{-0.2cm}
\end{figure}

\begin{figure}
    \centering
\includegraphics[scale=0.48]{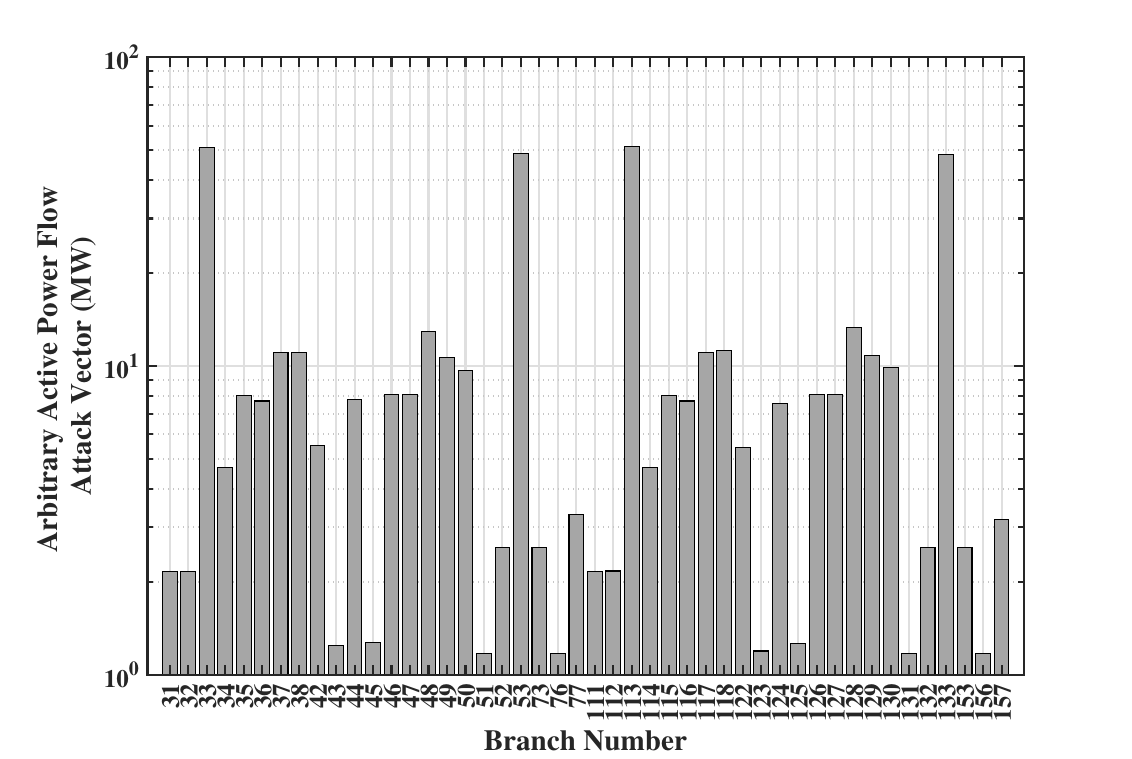}
	\caption{Active Power Flow attack vector for arbitrary scenario.}
	\label{fig:powerflow_arbitrary}
 \vspace{-.4cm}
\end{figure}

\begin{figure}
    \centering
\includegraphics[scale=0.48]{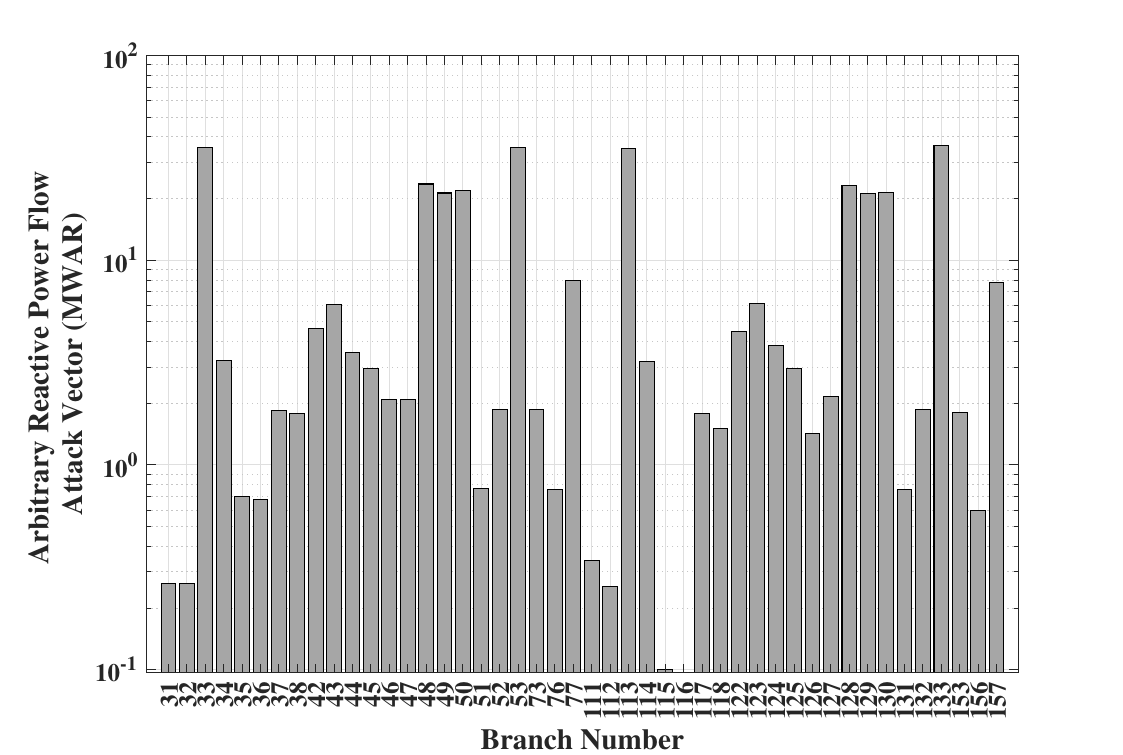}
	\caption{Reactive Power Flow attack vector for arbitrary scenario.}
	\label{fig:reactivepowerflow_arbitrary}
 \vspace{-.4cm}
\end{figure}

\mrn{Additionally, the active and reactive power injection attack vectors are depicted in Figs.~\ref{fig:powerinjection_arbitrary} and \ref{fig:reactivepowerinjection_arbitrary}. These figures show how the values of active and reactive power injection measurements for non-zero injection buses have been altered. It is important to note that the measurements for all other non-zero injection buses outside the attack zone should remain unchanged. Moreover, if all power flow constraints were truly modified, the measurements of zero injection buses within the attack zone should also remain unaffected.
Figs.~\ref{fig:powerflow_arbitrary} and \ref{fig:reactivepowerflow_arbitrary} illustrate the attack vectors for active and reactive power flow across all lines within the attack zone. It is evident that each line in the attack zone experiences a change in power flow due to variations in the voltage at one end (connected to boundary buses) or at both ends of the line.}

\begin{table}
\label{Table:arbitrary}
\caption{Corresponding values of Voltage magnitude and angle for before and after an arbitrary attack.}
\begin{tabular}{|c|cc|cc|}
\hline
                                                                                  & \multicolumn{2}{c|}{\textbf{\begin{tabular}[c]{@{}c@{}}Voltage \\ Magnitude\end{tabular}}}                                                                            & \multicolumn{2}{c|}{\textbf{\begin{tabular}[c]{@{}c@{}}Voltage \\ Angle\end{tabular}}}                                                                               \\ \cline{2-5} 
\multirow{-2}{*}{\textbf{\begin{tabular}[c]{@{}c@{}}Bus \\ Numbers\end{tabular}}} & \multicolumn{1}{c|}{\textbf{\begin{tabular}[c]{@{}c@{}}Before the   \\  Attack\end{tabular}}} & \textbf{\begin{tabular}[c]{@{}c@{}}Arbitrary \\  Attack\end{tabular}} & \multicolumn{1}{c|}{\textbf{\begin{tabular}[c]{@{}c@{}}Before \\  the  Attack\end{tabular}}} & \textbf{\begin{tabular}[c]{@{}c@{}}Arbitrary\\   Attack\end{tabular}} \\ \hline
\rowcolor[HTML]{C0C0C0} 
\textbf{20}                                                                       & \multicolumn{1}{c|}{\cellcolor[HTML]{C0C0C0}\textbf{0.9638}}                                  & \textbf{0.9638}                                                       & \multicolumn{1}{c|}{\cellcolor[HTML]{C0C0C0}\textbf{-0.2346}}                                & \textbf{-0.2346}                                                      \\ \hline
\rowcolor[HTML]{FFCCC9} 
\textbf{21}                                                                       & \multicolumn{1}{c|}{\cellcolor[HTML]{FFCCC9}\textbf{1.0085}}                                  & \textbf{1.0061}                                                       & \multicolumn{1}{c|}{\cellcolor[HTML]{FFCCC9}\textbf{-0.2257}}                                & \textbf{-0.2076}                                                      \\ \hline
\rowcolor[HTML]{FFCCC9} 
\textbf{22}                                                                       & \multicolumn{1}{c|}{\cellcolor[HTML]{FFCCC9}\textbf{1.0097}}                                  & \textbf{1.0086}                                                       & \multicolumn{1}{c|}{\cellcolor[HTML]{FFCCC9}\textbf{-0.2247}}                                & \textbf{-0.204}                                                       \\ \hline
\rowcolor[HTML]{FFCCC9} 
\textbf{23}                                                                       & \multicolumn{1}{c|}{\cellcolor[HTML]{FFCCC9}\textbf{1.0083}}                                  & \textbf{1.007}                                                        & \multicolumn{1}{c|}{\cellcolor[HTML]{FFCCC9}\textbf{-0.2258}}                                & \textbf{-0.194}                                                       \\ \hline
\rowcolor[HTML]{FFCCC9} 
\textbf{24}                                                                       & \multicolumn{1}{c|}{\cellcolor[HTML]{FFCCC9}\textbf{0.9992}}                                  & \textbf{0.9992}                                                       & \multicolumn{1}{c|}{\cellcolor[HTML]{FFCCC9}\textbf{-0.232}}                                 & \textbf{-0.232}                                                       \\ \hline
\rowcolor[HTML]{FFCCC9} 
\textbf{25}                                                                       & \multicolumn{1}{c|}{\cellcolor[HTML]{FFCCC9}\textbf{0.9825}}                                  & \textbf{1.0257}                                                       & \multicolumn{1}{c|}{\cellcolor[HTML]{FFCCC9}\textbf{-0.3172}}                                & \textbf{-0.3634}                                                      \\ \hline
\rowcolor[HTML]{FFCCC9} 
\textbf{26}                                                                       & \multicolumn{1}{c|}{\cellcolor[HTML]{FFCCC9}\textbf{0.9588}}                                  & \textbf{0.9735}                                                       & \multicolumn{1}{c|}{\cellcolor[HTML]{FFCCC9}\textbf{-0.2266}}                                & \textbf{-0.1938}                                                      \\ \hline
\rowcolor[HTML]{C0C0C0} 
\textbf{27}                                                                       & \multicolumn{1}{c|}{\cellcolor[HTML]{C0C0C0}\textbf{0.9815}}                                  & \textbf{0.9815}                                                       & \multicolumn{1}{c|}{\cellcolor[HTML]{C0C0C0}\textbf{-0.201}}                                 & \textbf{-0.201}                                                       \\ \hline
\rowcolor[HTML]{FFCCC9} 
\textbf{30}                                                                       & \multicolumn{1}{c|}{\cellcolor[HTML]{FFCCC9}\textbf{0.9627}}                                  & \textbf{1.0073}                                                       & \multicolumn{1}{c|}{\cellcolor[HTML]{FFCCC9}\textbf{-0.3267}}                                & \textbf{-0.3702}                                                      \\ \hline
\rowcolor[HTML]{FFCCC9} 
\textbf{31}                                                                       & \multicolumn{1}{c|}{\cellcolor[HTML]{FFCCC9}\textbf{0.9359}}                                  & \textbf{0.9998}                                                       & \multicolumn{1}{c|}{\cellcolor[HTML]{FFCCC9}\textbf{-0.3383}}                                & \textbf{-0.3477}                                                      \\ \hline
\rowcolor[HTML]{FFCCC9} 
\textbf{32}                                                                       & \multicolumn{1}{c|}{\cellcolor[HTML]{FFCCC9}\textbf{0.9499}}                                  & \textbf{1.0025}                                                       & \multicolumn{1}{c|}{\cellcolor[HTML]{FFCCC9}\textbf{-0.3231}}                                & \textbf{-0.3241}                                                      \\ \hline
\rowcolor[HTML]{FFCCC9} 
\textbf{33}                                                                       & \multicolumn{1}{c|}{\cellcolor[HTML]{FFCCC9}\textbf{0.9476}}                                  & \textbf{1.0018}                                                       & \multicolumn{1}{c|}{\cellcolor[HTML]{FFCCC9}\textbf{-0.3238}}                                & \textbf{-0.3235}                                                      \\ \hline
\rowcolor[HTML]{FFCCC9} 
\textbf{34}                                                                       & \multicolumn{1}{c|}{\cellcolor[HTML]{FFCCC9}\textbf{0.9592}}                                  & \textbf{0.9907}                                                       & \multicolumn{1}{c|}{\cellcolor[HTML]{FFCCC9}\textbf{-0.2469}}                                & \textbf{-0.2986}                                                      \\ \hline
\rowcolor[HTML]{FFCCC9} 
\textbf{35}                                                                       & \multicolumn{1}{c|}{\cellcolor[HTML]{FFCCC9}\textbf{0.9662}}                                  & \textbf{0.9932}                                                       & \multicolumn{1}{c|}{\cellcolor[HTML]{FFCCC9}\textbf{-0.2427}}                                & \textbf{-0.2972}                                                      \\ \hline
\rowcolor[HTML]{FFCCC9} 
\textbf{36}                                                                       & \multicolumn{1}{c|}{\cellcolor[HTML]{FFCCC9}\textbf{0.9758}}                                  & \textbf{0.996}                                                        & \multicolumn{1}{c|}{\cellcolor[HTML]{FFCCC9}\textbf{-0.238}}                                 & \textbf{-0.2741}                                                      \\ \hline
\rowcolor[HTML]{FFCCC9} 
\textbf{37}                                                                       & \multicolumn{1}{c|}{\cellcolor[HTML]{FFCCC9}\textbf{0.9849}}                                  & \textbf{1.0004}                                                       & \multicolumn{1}{c|}{\cellcolor[HTML]{FFCCC9}\textbf{-0.2347}}                                & \textbf{-0.26}                                                        \\ \hline
\rowcolor[HTML]{C0C0C0} 
\textbf{38}                                                                       & \multicolumn{1}{c|}{\cellcolor[HTML]{C0C0C0}\textbf{1.0128}}                                  & \textbf{1.0128}                                                       & \multicolumn{1}{c|}{\cellcolor[HTML]{C0C0C0}\textbf{-0.2223}}                                & \textbf{-0.2223}                                                      \\ \hline
\rowcolor[HTML]{FFCCC9} 
\textbf{39}                                                                       & \multicolumn{1}{c|}{\cellcolor[HTML]{FFCCC9}\textbf{0.9828}}                                  & \textbf{0.9989}                                                       & \multicolumn{1}{c|}{\cellcolor[HTML]{FFCCC9}\textbf{-0.2355}}                                & \textbf{-0.2605}                                                      \\ \hline
\rowcolor[HTML]{FFCCC9} 
\textbf{40}                                                                       & \multicolumn{1}{c|}{\cellcolor[HTML]{FFCCC9}\textbf{0.9728}}                                  & \textbf{0.993}                                                        & \multicolumn{1}{c|}{\cellcolor[HTML]{FFCCC9}\textbf{-0.2384}}                                & \textbf{-0.2726}                                                      \\ \hline
\rowcolor[HTML]{C0C0C0} 
\textbf{56}                                                                       & \multicolumn{1}{c|}{\cellcolor[HTML]{C0C0C0}\textbf{0.9684}}                                  & \textbf{0.9684}                                                       & \multicolumn{1}{c|}{\cellcolor[HTML]{C0C0C0}\textbf{-0.2804}}                                & \textbf{-0.2804}                                                      \\ \hline
\rowcolor[HTML]{FFCCC9} 
\textbf{57}                                                                       & \multicolumn{1}{c|}{\cellcolor[HTML]{FFCCC9}\textbf{0.9648}}                                  & \textbf{0.9914}                                                       & \multicolumn{1}{c|}{\cellcolor[HTML]{FFCCC9}\textbf{-0.2894}}                                & \textbf{-0.2943}                                                      \\ \hline
\end{tabular}
\end{table}

\textbf{Sparse Attack:}
\mrn{In this scenario, we again consider buses 20, 21, 22, 23, 24, 25, 26, 27, 30, 31, 32, 33, 34, 35, 36, 37, 38, 39, 40, 56 and 57 as part of the attack zone. Prior to solving the problem, buses 20, 27, 38 and 56 are designated as boundary buses. To control the power flow within the attack zone, the voltages at these boundary buses must be fixed. Similar to the previous scenario we enforce Equation~\ref{eq:active_overload} with $W=3$.}

\mrn{In this scenario, we design a sparse AC FDI attack by solving Equations~\eqref{eq:obj1}-\eqref{eq:reactive_overload}. In these constraints, when the variable $z$ is equal to zero, the corresponding buses are fixed at a predetermined voltage, which updates the attack zone and boundary buses accordingly. As a result, when the solver processes the problem, buses 
21, 22, 23 and 25 are identified as buses that can be altered.
As shown in Table \ref{Table:sparse}, the number of selected buses that need to be altered to execute a successful attack is reduced compared to the arbitrary scenario, thus achieving a sparse AC FDI attack. In this table, gray rows represent boundary buses with fixed voltages, white rows are flexible nodes in the attack zone that could be selected for attack but remain fixed after solving the problem, and pink rows are flexible nodes in the attack zone that have been selected for manipulation after solving the problem. Note that the proposed algorithm achieves a successful AC FDI attack with the same impact on the grid by manipulating significantly fewer measurements.}

\begin{small}
\vspace{-.1cm}
\begin{table}
\caption{Corresponding values of Voltage magnitude and angle for before and after sparse attack.}
\label{Table:sparse}
\begin{tabular}{|c|cc|cc|}
\hline
                                                                                  & \multicolumn{2}{c|}{\textbf{\begin{tabular}[c]{@{}c@{}}Voltage \\ Magnitude\end{tabular}}}                                                                            & \multicolumn{2}{c|}{\textbf{\begin{tabular}[c]{@{}c@{}}Voltage \\ Angle\end{tabular}}}                                                                               \\ \cline{2-5} 
\multirow{-2}{*}{\textbf{\begin{tabular}[c]{@{}c@{}}Bus \\ Numbers\end{tabular}}} & \multicolumn{1}{c|}{\textbf{\begin{tabular}[c]{@{}c@{}}Before the   \\  Attack\end{tabular}}} & \textbf{\begin{tabular}[c]{@{}c@{}}Arbitrary \\  Attack\end{tabular}} & \multicolumn{1}{c|}{\textbf{\begin{tabular}[c]{@{}c@{}}Before \\  the  Attack\end{tabular}}} & \textbf{\begin{tabular}[c]{@{}c@{}}Arbitrary\\   Attack\end{tabular}} \\ \hline
\rowcolor[HTML]{C0C0C0} 
\textbf{20}                                                                       & \multicolumn{1}{c|}{\cellcolor[HTML]{C0C0C0}\textbf{0.9638}}                                  & \textbf{0.9638}                                                       & \multicolumn{1}{c|}{\cellcolor[HTML]{C0C0C0}\textbf{-0.2346}}                                & \textbf{-0.2346}                                                      \\ \hline
\rowcolor[HTML]{FFCCC9} 
\textbf{21}                                                                       & \multicolumn{1}{c|}{\cellcolor[HTML]{FFCCC9}\textbf{1.0085}}                                  & \textbf{1.0061}                                                       & \multicolumn{1}{c|}{\cellcolor[HTML]{FFCCC9}\textbf{-0.2257}}                                & \textbf{-0.2076}                                                      \\ \hline
\rowcolor[HTML]{FFCCC9} 
\textbf{22}                                                                       & \multicolumn{1}{c|}{\cellcolor[HTML]{FFCCC9}\textbf{1.0097}}                                  & \textbf{1.0086}                                                       & \multicolumn{1}{c|}{\cellcolor[HTML]{FFCCC9}\textbf{-0.2247}}                                & \textbf{-0.204}                                                       \\ \hline
\rowcolor[HTML]{FFCCC9} 
\textbf{23}                                                                       & \multicolumn{1}{c|}{\cellcolor[HTML]{FFCCC9}\textbf{1.0083}}                                  & \textbf{1.007}                                                        & \multicolumn{1}{c|}{\cellcolor[HTML]{FFCCC9}\textbf{-0.2258}}                                & \textbf{-0.194}                                                       \\ \hline
\rowcolor[HTML]{FFFFFF} 
\textbf{24}                                                                       & \multicolumn{1}{c|}{\cellcolor[HTML]{FFFFFF}\textbf{0.9992}}                                  & \textbf{0.9992}                                                       & \multicolumn{1}{c|}{\cellcolor[HTML]{FFFFFF}\textbf{-0.232}}                                 & \textbf{-0.232}                                                       \\ \hline
\rowcolor[HTML]{FFCCC9} 
\textbf{25}                                                                       & \multicolumn{1}{c|}{\cellcolor[HTML]{FFCCC9}\textbf{0.9825}}                                  & \textbf{1.0257}                                                       & \multicolumn{1}{c|}{\cellcolor[HTML]{FFCCC9}\textbf{-0.3172}}                                & \textbf{-0.3634}                                                      \\ \hline
\rowcolor[HTML]{FFFFFF} 
\textbf{26}                                                                       & \multicolumn{1}{c|}{\cellcolor[HTML]{FFFFFF}\textbf{0.9588}}                                  & \cellcolor[HTML]{FFFFFF}\textbf{0.9588}                               & \multicolumn{1}{c|}{\cellcolor[HTML]{FFFFFF}\textbf{-0.2266}}                                & \textbf{-0.2266}                                                      \\ \hline
\rowcolor[HTML]{C0C0C0} 
\textbf{27}                                                                       & \multicolumn{1}{c|}{\cellcolor[HTML]{C0C0C0}\textbf{0.9815}}                                  & \textbf{0.9815}                                                       & \multicolumn{1}{c|}{\cellcolor[HTML]{C0C0C0}\textbf{-0.201}}                                 & \textbf{-0.201}                                                       \\ \hline
\rowcolor[HTML]{FFFFFF} 
\textbf{30}                                                                       & \multicolumn{1}{c|}{\cellcolor[HTML]{FFFFFF}\textbf{0.9627}}                                  & \cellcolor[HTML]{FFFFFF}\textbf{0.9627}                               & \multicolumn{1}{c|}{\cellcolor[HTML]{FFFFFF}\textbf{-0.3267}}                                & \cellcolor[HTML]{FFFFFF}\textbf{-0.3267}                              \\ \hline
\rowcolor[HTML]{FFFFFF} 
\textbf{31}                                                                       & \multicolumn{1}{c|}{\cellcolor[HTML]{FFFFFF}\textbf{0.9359}}                                  & \cellcolor[HTML]{FFFFFF}\textbf{0.9359}                               & \multicolumn{1}{c|}{\cellcolor[HTML]{FFFFFF}\textbf{-0.3383}}                                & \cellcolor[HTML]{FFFFFF}\textbf{-0.3383}                              \\ \hline
\rowcolor[HTML]{FFFFFF} 
\textbf{32}                                                                       & \multicolumn{1}{c|}{\cellcolor[HTML]{FFFFFF}\textbf{0.9499}}                                  & \cellcolor[HTML]{FFFFFF}\textbf{0.9499}                               & \multicolumn{1}{c|}{\cellcolor[HTML]{FFFFFF}\textbf{-0.3231}}                                & \cellcolor[HTML]{FFFFFF}\textbf{-0.3231}                              \\ \hline
\rowcolor[HTML]{FFFFFF} 
\textbf{33}                                                                       & \multicolumn{1}{c|}{\cellcolor[HTML]{FFFFFF}\textbf{0.9476}}                                  & \cellcolor[HTML]{FFFFFF}\textbf{0.9476}                               & \multicolumn{1}{c|}{\cellcolor[HTML]{FFFFFF}\textbf{-0.3238}}                                & \cellcolor[HTML]{FFFFFF}\textbf{-0.3238}                              \\ \hline
\rowcolor[HTML]{FFFFFF} 
\textbf{34}                                                                       & \multicolumn{1}{c|}{\cellcolor[HTML]{FFFFFF}\textbf{0.9592}}                                  & \cellcolor[HTML]{FFFFFF}\textbf{0.9592}                               & \multicolumn{1}{c|}{\cellcolor[HTML]{FFFFFF}\textbf{-0.2469}}                                & \cellcolor[HTML]{FFFFFF}\textbf{-0.2469}                              \\ \hline
\rowcolor[HTML]{FFFFFF} 
\textbf{35}                                                                       & \multicolumn{1}{c|}{\cellcolor[HTML]{FFFFFF}\textbf{0.9662}}                                  & \cellcolor[HTML]{FFFFFF}\textbf{0.9662}                               & \multicolumn{1}{c|}{\cellcolor[HTML]{FFFFFF}\textbf{-0.2427}}                                & \cellcolor[HTML]{FFFFFF}\textbf{-0.2427}                              \\ \hline
\rowcolor[HTML]{FFFFFF} 
\textbf{36}                                                                       & \multicolumn{1}{c|}{\cellcolor[HTML]{FFFFFF}\textbf{0.9758}}                                  & \cellcolor[HTML]{FFFFFF}\textbf{0.9758}                               & \multicolumn{1}{c|}{\cellcolor[HTML]{FFFFFF}\textbf{-0.238}}                                 & \cellcolor[HTML]{FFFFFF}\textbf{-0.238}                               \\ \hline
\rowcolor[HTML]{FFFFFF} 
\textbf{37}                                                                       & \multicolumn{1}{c|}{\cellcolor[HTML]{FFFFFF}\textbf{0.9849}}                                  & \cellcolor[HTML]{FFFFFF}\textbf{0.9849}                               & \multicolumn{1}{c|}{\cellcolor[HTML]{FFFFFF}\textbf{-0.2347}}                                & \cellcolor[HTML]{FFFFFF}\textbf{-0.2347}                              \\ \hline
\rowcolor[HTML]{C0C0C0} 
\textbf{38}                                                                       & \multicolumn{1}{c|}{\cellcolor[HTML]{C0C0C0}\textbf{1.0128}}                                  & \textbf{1.0128}                                                       & \multicolumn{1}{c|}{\cellcolor[HTML]{C0C0C0}\textbf{-0.2223}}                                & \textbf{-0.2223}                                                      \\ \hline
\rowcolor[HTML]{FFFFFF} 
\textbf{39}                                                                       & \multicolumn{1}{c|}{\cellcolor[HTML]{FFFFFF}\textbf{0.9828}}                                  & \cellcolor[HTML]{FFFFFF}\textbf{0.9828}                               & \multicolumn{1}{c|}{\cellcolor[HTML]{FFFFFF}\textbf{-0.2355}}                                & \cellcolor[HTML]{FFFFFF}\textbf{-0.2355}                              \\ \hline
\rowcolor[HTML]{FFFFFF} 
\textbf{40}                                                                       & \multicolumn{1}{c|}{\cellcolor[HTML]{FFFFFF}\textbf{0.9728}}                                  & \cellcolor[HTML]{FFFFFF}\textbf{0.9728}                               & \multicolumn{1}{c|}{\cellcolor[HTML]{FFFFFF}\textbf{-0.2384}}                                & \cellcolor[HTML]{FFFFFF}\textbf{-0.2384}                              \\ \hline
\rowcolor[HTML]{C0C0C0} 
\textbf{56}                                                                       & \multicolumn{1}{c|}{\cellcolor[HTML]{C0C0C0}\textbf{0.9684}}                                  & \textbf{0.9684}                                                       & \multicolumn{1}{c|}{\cellcolor[HTML]{C0C0C0}\textbf{-0.2804}}                                & \textbf{-0.2804}                                                      \\ \hline
\rowcolor[HTML]{FFFFFF} 
\textbf{57}                                                                       & \multicolumn{1}{c|}{\cellcolor[HTML]{FFFFFF}\textbf{0.9648}}                                  & \cellcolor[HTML]{FFFFFF}\textbf{0.9648}                               & \multicolumn{1}{c|}{\cellcolor[HTML]{FFFFFF}\textbf{-0.2894}}                                & \cellcolor[HTML]{FFFFFF}\textbf{-0.2894}                              \\ \hline
\end{tabular}
\end{table}
\end{small}

\begin{figure}
% \vspace{-1.5cm}
    \centering
\includegraphics[scale=0.47]{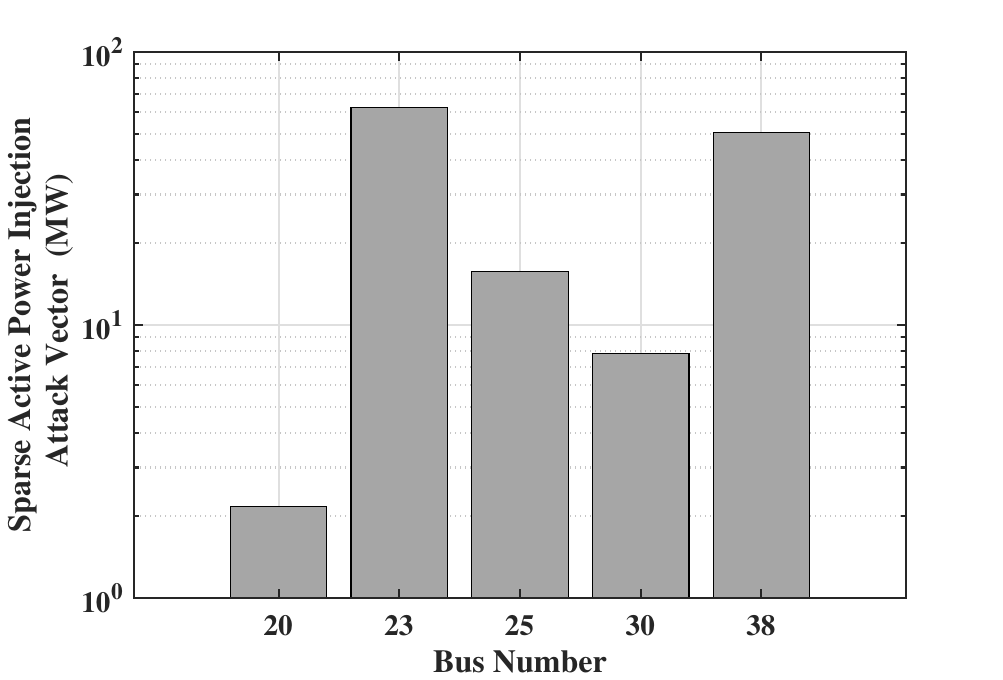}
	\caption{\mrn{Active power injection attack vector in the sparse AC FDI attack scenario. Bus numbers indicate the corresponding measurement numbers on the buses. Compared to Fig.~\ref{fig:powerinjection_arbitrary}, it is clear that the number of manipulated measurements in the sparse attack scenario is fewer than in the arbitrary FDI attack. }}
	\label{fig:active_injection_sparse}
\vspace{0cm}
\end{figure}

\begin{figure}
 \vspace{-0.0cm}
    \centering
\includegraphics[scale=0.47]{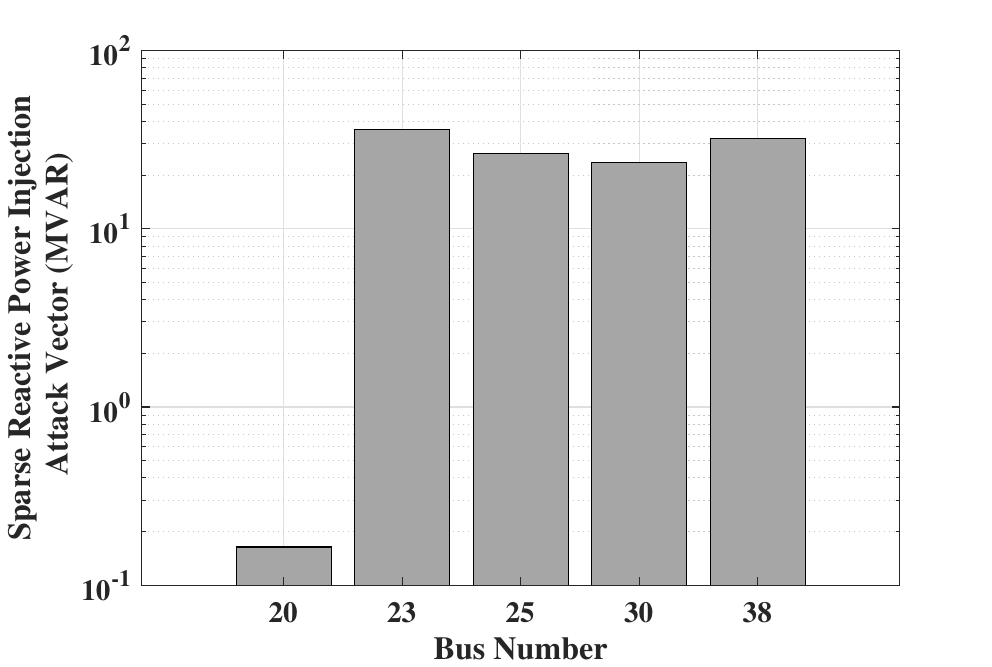}
	\caption{\mrn{Reactive power injection attack vector in the sparse AC FDI attack scenario. Bus numbers indicate the corresponding measurement numbers on the buses.  Compared to Fig.~\ref{fig:reactivepowerinjection_arbitrary}, it is clear that the number of manipulated measurements in the sparse attack scenario is fewer than in the arbitrary FDI attack.}}
\label{fig:reactive_injection_sparse}
%\vspace{-1cm}
\end{figure}

\begin{figure}
    \centering
\includegraphics[scale=0.47]{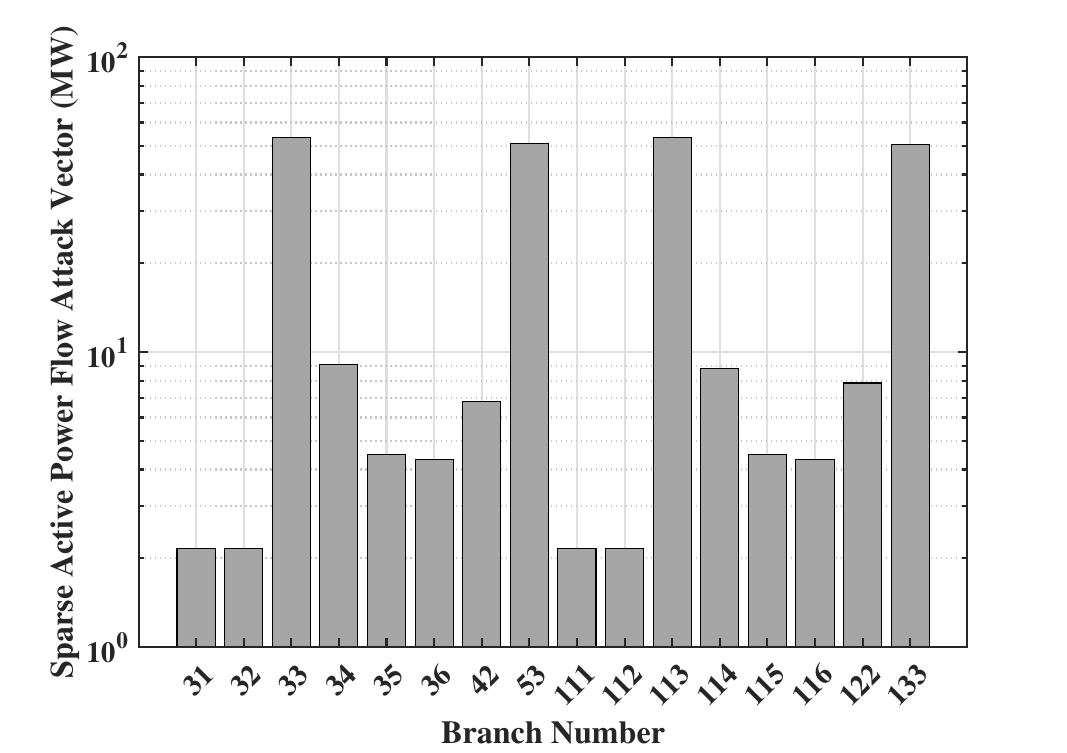}
	\caption{\mrn{Active power flow attack vector associated with the sparse AC FDI attack scenario. Branch numbers indicate the corresponding measurement numbers on the lines in the attack zone. }}
\label{fig:active_flow_sparse}
%\vspace{-.1cm}
\end{figure}

\begin{figure}
    \centering
\includegraphics[scale=0.47]{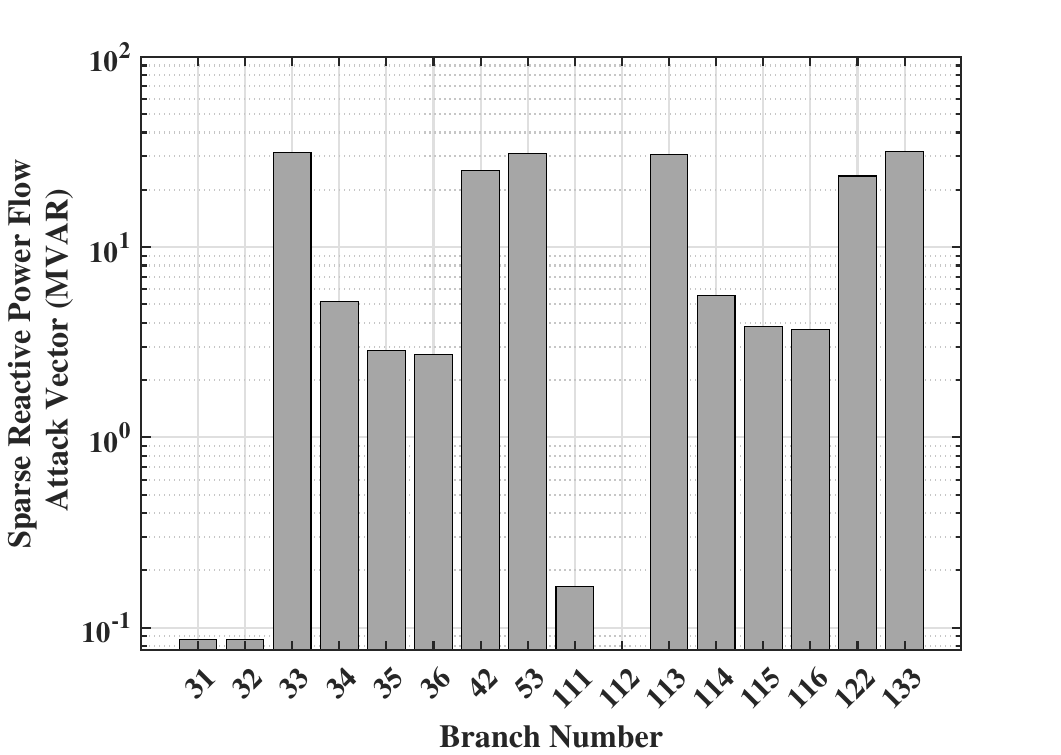}
	\caption{\mrn{Reactive Power flow attack vector in sparse AC FDI attack scenario. Branch numbers indicate the corresponding measurements' number on the lines in the attack zone. Compared to Fig.~\ref{fig:reactivepowerflow_arbitrary}, it is clear that the number of manipulated measurements in the sparse attack scenario is fewer than in the arbitrary FDI attack.} }
\label{fig:reactive_flow_sparse}
\end{figure}

\mrn{Figs.~\ref{fig:active_injection_sparse} and \ref{fig:reactive_injection_sparse} display the active and reactive power injections at non-zero injection buses in the new attack zone, while Figs.~\ref{fig:active_flow_sparse} and \ref{fig:reactive_flow_sparse} show the active and reactive power flow in the lines within the sparse attack zone.}

\section{Conclusion}
\label{sec:conclusion}
\mrn{This paper addresses the design and implementation of sparse AC FDI attacks, diverging from previous studies that predominantly focused on DC FDI attacks. This work tackles the complexity and nonlinearity inherent in AC power flows, which presents unique challenges for FDI attack studies.
We formulated the problem as a Mixed Integer Nonlinear Programming (MINLP) problem, aiming to strategically select a minimal set of measurements for manipulation while preserving the non-linearity of AC power flow equations. We employed binary variables to indicate the selection of measurements and continuous variables to represent their values, allowing dynamic adjustment of power flow constraints based on the manipulation status. The big-M method and conditional constraints were effectively used to manage both fixed and variable measurement parameters. Specifically, we fixed the voltage of unselected measurements and allowed variable values for selected PMU measurements. By integrating these variables into a mixed integer nonlinear optimization problem, we minimized the number of altered measurements, achieving a sparse AC FDI attack.
Simulation results on the IEEE 57-bus test system demonstrate the effectiveness of the proposed approach. The results reveal that a minimal number of measurements can be manipulated to significantly impact state estimation, underscoring the potential vulnerability of power systems to such attacks. However, solving the MINLP presents a challenge, as the increased complexity can affect efficiency. Future work will focus on relaxing the nonlinear power flow equations to improve the efficacy of the proposed method.}

\bibliographystyle{IEEEtran}
\IEEEtriggeratref{40}
\bibliography{ref}
\end{document}